\documentclass[letterpaper, 10pt, conference]{ieeeconf}

\IEEEoverridecommandlockouts

\let\labelindent\relax
\usepackage{cite}
\usepackage{amsmath,amssymb,amsfonts}
\usepackage{graphicx}
\usepackage{enumerate}
\usepackage{enumitem}
\usepackage{textcomp}
\usepackage{xcolor}
\def\BibTeX{{\rm B\kern-.05em{\sc i\kern-.025em b}\kern-.08em
    T\kern-.1667em\lower.7ex\hbox{E}\kern-.125emX}}

\usepackage{newtxmath, newtxtext}

\usepackage{amsmath}
\usepackage{amssymb}
\usepackage{mathtools}
\usepackage{amsthm}
\usepackage{ragged2e}
\usepackage{afterpage}
\usepackage{float}
\usepackage{algorithm}
\usepackage{algpseudocode}
\usepackage{placeins}
\usepackage{mathrsfs}
\usepackage{booktabs}
\usepackage{balance}
\usepackage{dutchcal}

\theoremstyle{plain}
\newtheorem{theorem}{Theorem}[section]
\newtheorem{proposition}[theorem]{Proposition}
\newtheorem{lemma}[theorem]{Lemma}

\theoremstyle{definition}
\newtheorem{definition}[theorem]{Definition}

\newtheorem{remark}[theorem]{Remark}

\newcommand{\norm}[1]{\left\lVert#1\right\rVert}

\title{\LARGE \bf 
Control of Markov Jump Linear Systems\\with Uncertain Lumpable Cluster Observations
}

\author{Carlos A. F. Persiani, Ram Padmanabhan, Melkior Ornik and Marco H. Terra 
\thanks{Carlos A. F. Persiani and Marco H. Terra are with the University of São Paulo at São Carlos, SP, Brazil. Emails: \texttt{carlospersiani@usp.br, terra@sc.usp.br}. Their work was supported by the Coordenação de Aperfeiçoamento de Pessoal de Nível Superior - Brasil (CAPES) --- Finance Code 001. Corresponding Author: Carlos A. F. Persiani.
}
\thanks{Ram Padmanabhan and Melkior Ornik are with the University of Illinois Urbana-Champaign, Urbana, IL 61801, USA. Emails: \texttt{\{ramp3, mornik\}@illinois.edu}. Their work was supported by the Air Force Office of Scientific Research under Grant FA9550-23-1-0131.}
\thanks{This work was conducted while Carlos A. F. Persiani was a visiting research scholar at University of Illinois Urbana-Champaign, Urbana, IL 61801, USA.}
}

\begin{document}

\maketitle

\begin{abstract}
In this paper, we consider the control of Markov jump linear systems when the active system mode is not exactly known. We rather assume that this true mode is only known to belong to an observed cluster of modes, making both the dynamics and the associated transition probabilities uncertain. We construct a min-max optimal control problem in this setting, with the objective of robust regulation. Using results in regularized least mean-square optimization, we describe how the solution of this problem can be reduced to a Riccati recursion, such that the controller depends only on knowledge of which cluster the mode lies in. We subsequently establish convergence of the recursion and stability of the resulting closed-loop system. A numerical example illustrates the effectiveness of the proposed methodology.
\end{abstract}

\section{Introduction}\label{section:int}

Discrete-time Markov jump linear systems constitute an important class of stochastic hybrid systems in which the switching mechanism between modes is governed by a Markov chain. Such models arise naturally in systems subject to failures, changes in operating conditions, network induced effects, and other random switching phenomena. Owing to their relevance, fundamental problems concerning stability, stabilization, optimal and robust control have been extensively investigated \cite{costa2005discrete}. A substantial portion of this theory, however, relies on the assumption that the true mode is perfectly available to the controller.

In many applications, complete knowledge of the mode may not be possible, and only information about a cluster of candidate modes can be obtained. An important formulation for this setting was introduced in \cite{do2002h2}, where the modes are partitioned into clusters and the controller has access only to the cluster containing the true mode. A controller is assigned to each cluster of modes, and stabilizing $\mathcal{H}_2$ controllers are obtained through linear matrix inequalities. Related formulations have subsequently considered detector-based observations \cite{do2014detector}, $\mathcal{H}_\infty$ performance under partial observations \cite{todorov2018detector}, unknown transition probabilities \cite{li2019passivity, zhu2020asynchronous} and modes observed at random times \cite{7170988, ogura2018state}.

The property of \emph{lumpability} provides conditions under which some clusters can still be represented by a Markov chain \cite{jernigan2003testing, peng1996weak}. In particular, under \emph{strong lumpability}, modes belonging to the same cluster have identical transition probabilities toward every other cluster, yielding a well-defined transition probability between clusters. The property of \emph{quasi-lumpability} relaxes this requirement by allowing these probabilities to differ within a prescribed tolerance (see \cite{marin2022proportional} and references therein). When lumpability does not hold, the transition probability from one cluster to another is no longer uniquely determined by the current cluster alone, but may depend on the particular mode occupied within it. From the perspective of an observer with access only to the clustered mode, this dependence introduces uncertainty into the clustered transition probabilities, such that the clustered dynamics can be represented by an uncertain transition matrix, as illustrated in the example below.

\begin{figure}[!h]
    \centering
    \includegraphics[width=0.9\linewidth]{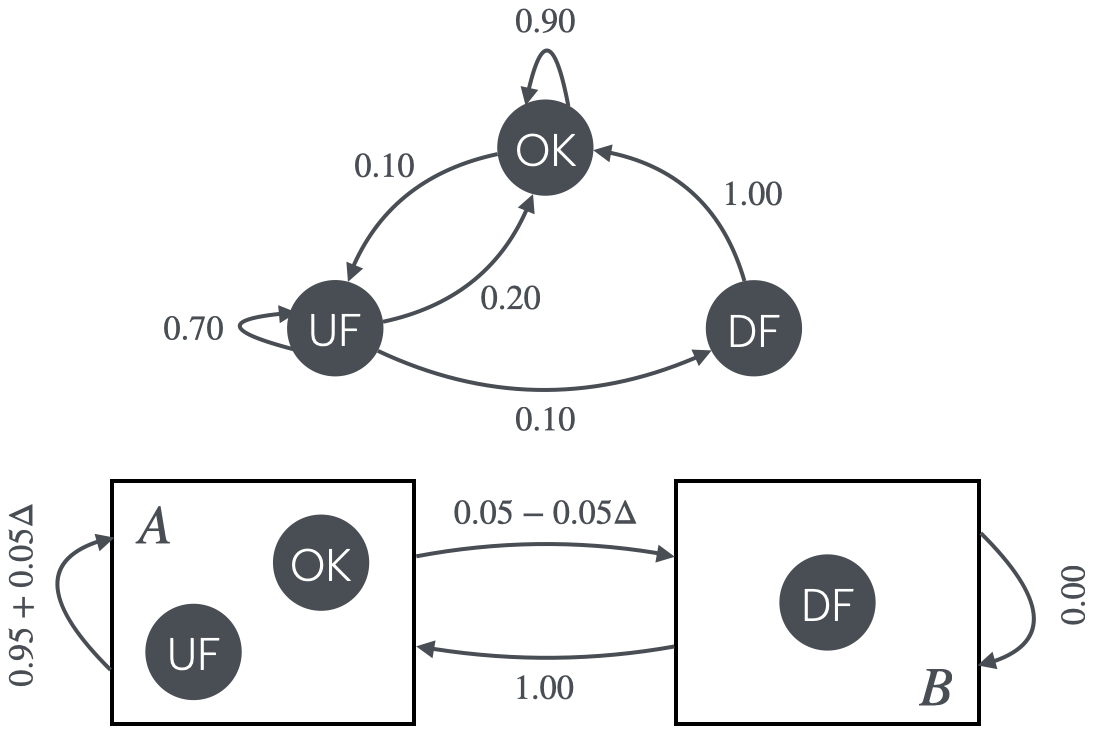}
    \caption{The transition between clusters is parameterized by an uncertain quantity $\Delta\in \{-1,1\}$. For $\Delta=1$ ($\Delta=-1$), cluster $A$ coincides with mode OK (UF). In this particular example, cluster $B$ has no associated uncertainty.}
    \label{fig::MarkovClusterUncertain}
    \vspace{-0.2cm}
\end{figure}

\textbf{Example:} Consider Example 6.1 from \cite{todorov2018detector} of a Markov jump linear system with three modes: Normal Operation (OK), Undetected Fault (UF) and Detected Fault (DF). These modes are \emph{aggregated} into two clusters: Apparent Normal (A), containing modes OK and UF, and Apparent Fault (B), containing mode DF. Figure \ref{fig::MarkovClusterUncertain} illustrates both the original transition probability matrix and its clustered uncertain counterpart.

Existing approaches to the control of Markov jump linear systems with cluster observations do not account for the uncertainty that naturally arises in the stochastic transition matrix of the clusters. This paper seeks to address this issue, where we model the transition probabilities between clusters through an uncertain stochastic matrix similar to Example 1. There thus exist two distinct sources of uncertainty, in that the true mode is unknown and only the cluster it belongs to is known, and the transition probabilities between clusters are themselves uncertain. In this setting, we formulate a robust optimal control problem, where the uncertain effects seek to maximize a quadratic objective while the control input aims to minimize this objective. {Building upon results} from \cite{cerri2017recursive, bueno2022regulation, odorico2025robust}, we show that the solution can be recast as a set of coupled Riccati equations, from which we derive a cluster-dependent control law, providing a closed-form solution in contrast to existing approaches in literature. We subsequently provide an analysis of stability and convergence for these coupled Riccati recursions, and demonstrate the effectiveness of this approach through a simulation example.



The remainder of this paper is organized as follows: Section \ref{sec::problem_stat} presents the problem statement and formulates the robust optimal control problem we aim to solve. Section \ref{sec::preliminaries} introduces a modified version of the central lemma in \cite{cerri2017recursive} for solving a particular min-max quadratic program.  Section \ref{sec::RecursiveRegulation} shows relaxations of the original problem and, using results from Section \ref{sec::preliminaries}, derives a closed-form solution to the resulting relaxed problem. Section \ref{sec::Stability} provides conditions under which stability and convergence properties hold, and Section \ref{sec::Examples} presents a numerical example illustrating the performance of the proposed method.

\textbf{Notation and Properties}: The $n$-dimensional identity matrix is denoted $I_n$. The expectation of a random variable is denoted $\mathbb{E}\{\cdot\}$. If a matrix $P$ is symmetric positive (semi)definite, we write $P \succ 0$ ($P \succeq 0$). The Kronecker product between two matrices $A$ and $B$ is denoted $A\otimes B$. The general Euclidean $2$-norm is denoted by $\|\cdot\|$. The weighted squared Euclidean norm is denoted by $\norm{x}_{Q}^2 = x_{k}^{T}Q x_{k}$. The quantities $\rho_{\min}(W)$ and $\rho_{\max}(W)$ are the minimum and maximum eigenvalues {in norm} of a square matrix $W$. We recall the following properties from \cite{horn2012matrix}:

(F1) Woodbury Matrix Identity: For appropriate dimensions of matrices $A, X, R, Y$, $\left( A + XRY \right)^{-1}=A^{-1}-A^{-1}X \left( R^{-1} + Y A^{-1}X\right)^{-1} Y A^{-1}$ when the inverses are defined.

(F2) For a square matrix $A$ and $(I+A)^{-1}$ defined, $I - (I + A)^{-1} = A(I+A)^{-1}$.

(F3) Rayleigh Inequality: For $W \succeq 0$, $\rho_{\min}(W)\norm{x}^2 \leq \norm{x}^2_{W} \leq \rho_{\max}(W)\norm{x}^2$.

(F4) We say that a square matrix $\Delta$ is a \emph{contraction} if $\rho_{max}(\Delta)\leq1$. {A matrix $\Delta$ is a contraction if and only if $\Delta^T\Delta\preceq I$.}
\section{Problem Statement}\label{sec::problem_stat}

\subsection{System Dynamics and Clustering}
Consider the discrete-time Markov jump linear system \cite{costa2005discrete}
\begin{equation}\label{eq::system_dynamics}
    x_{k+1} = A_{\theta_{k}} x_{k} + B_{\theta_{k}}  u_{k},
\end{equation}
where $x_{k} \in \mathbb{R}^{n}$ is the state and $u_{k} \in \mathbb{R}^{m}$ is the control input. 
Let $\left\{ \theta_{k}\right\}_{k\geq0}$ be a Markov chain where each $\theta_k$ takes values in the set of modes $\Theta = \left\{ 1, ..., s\right\}$. The transition between modes of this Markov chain is governed by the nominal transition probability matrix $\mathbb{P} = \left[ p_{ij} \right] \in \mathbb{R}^{s \times s}$ such that
\begin{equation}\label{eq::prob_Markov_01}
    \mathbb{P}(\theta_{k+1} = j| \theta_{k} = i) = p_{ij}, \hspace{0.5cm} \mathbb{P}(\theta_{0}=i) = p_{i0},
\end{equation}
\begin{equation}\label{eq::prob_Markov_02}
    \sum_{j \in \Theta} p_{ij} = 1, \hspace{0.5cm} 0 \leq p_{ij} \leq 1.
\end{equation}

Let $\hat{\Theta}=\{1,\ldots,h\}$ denote the index set of observable clusters. Associated with each index $\hat{\theta}_k\in\hat{\Theta}$ is a cluster $L_{\hat{\theta}_k}\subseteq\Theta$, such that the collection of clusters $\mathcal{L}=\{L_1,\ldots,L_h\}$ satisfies $\bigcup_{l\in \hat{\Theta}}L_l=\Theta$. The process corresponding to clusters is denoted by $\{\hat{\theta}_k\}_{k\geq0}$, and it is assumed that the true Markov mode belongs to an observed cluster, i.e., $\theta_k\in L_{\hat{\theta}_k}$. We assume the transition between clusters is governed by an uncertain transition probability matrix $\mathbb{L}=[q_{ld}^\delta]\in \mathbb{R}^{h\times h}$, such that
\begin{equation}\label{eq::prob_Markov_cluster_01}
    \mathbb{L}(\hat{\theta}_{k+1} = d| \hat{\theta}_{k} = l) = q_{ld}^{\delta}, \hspace{0.5cm} \mathbb{L}(\hat{\theta}_{0}=l) = q_{l0},
\end{equation}
\begin{equation}\label{eq::prob_Markov_cluster_02}
    \sum_{d \in \hat{\Theta}} q_{ld}^\delta = 1, \hspace{0.5cm} 0 \leq q_{ld}^\delta \leq 1.
\end{equation}
Here, the superscript $\delta$ simply denotes uncertainty in the transition $q_{ld}$. Note that for a given observed cluster $l = \hat{\theta}_k$, the transition probability toward cluster $\hat{\theta}_{k+1} = d$ depends on the \emph{unknown} true mode $i\in L_l$. This transition probability, denoted $q_{ld}^{(i)}$, is given by
\begin{equation}
    q_{ld}^{(i)}=\sum_{j\in L_d} p_{ij},
\end{equation}
which corresponds to the sum of the transition probabilities from mode $i$ to all modes belonging to cluster $L_d$. Since the controller observes only cluster $L_l$, and not a particular mode $i$ within it, the corresponding cluster transition probability $q_{ld}^{\delta}$ is uncertain, not uniquely determined, and satisfies
\begin{equation}
    q_{ld}^{\delta}\in\left\{q_{ld}^{(i)}:i\in L_l\right\}.
\end{equation}

\begin{remark}
    If we assume the property of \emph{strong lumpability} \cite{jernigan2003testing} holds for the original Markov chain $\{\theta_k\}_{k \geq 0}$ with respect to an observed cluster partition, the transition probability among clusters $\mathbb{L}$ is uniquely determined and therefore does not exhibit uncertainty. In other words, $q_{ld}^{(i)}$ is the same quantity for all $i\in L_{l}$ for all $l\in \hat{\Theta}$. However, we do not make this assumption in our work. 
\end{remark}

Other works such as \cite{todorov2018detector} require knowledge of \emph{emission probability}, i.e., the distribution of clusters given a specific current mode. However, we do not assume this knowledge either.

\subsection{Control Objective}
Under the setup discussed above, we now formulate the robust optimal control problem that is the focus of this paper. We remark that as future system modes are unknown, we penalize the state evolution associated with each possible future mode by \emph{lifting} the system dynamics. 

Let $e_i\in\mathbb{R}^{s}$ denote the $i$-th canonical basis vector, and define
\begin{equation}\label{eq::lifted_matrices}
    F = \sum_{i=1}^{s} e_i \otimes A_{i}; \hspace{0.5em}
    G = \sum_{i=1}^{s} e_i \otimes B_{i}; \hspace{0.5em}
    C_i = e_i^{T}\otimes I_n.
\end{equation}
Using these matrices, the system dynamics \eqref{eq::system_dynamics} can be equivalently expressed in the lifted form
\begin{equation}\label{eq::open_dynamics}
\begin{aligned}
    z_{k+1} &= F x_k + G u_k,\\
    x_{k+1} &= C_{\hat{\theta}_k}^{\delta} z_{k+1},
\end{aligned}
\end{equation}
where $z_{k+1}\in\mathbb{R}^{sn}$ is the \emph{lifted state} collecting the evolution associated with all $s$ true modes. The matrix $C_{\hat{\theta}_k}^{\delta}$ is an uncertain output matrix belonging to the set $\left\{C_i:i\in L_{\hat{\theta}_k}\right\}$, where $C_i$ as defined in \eqref{eq::lifted_matrices} seeks to extract the state evolution associated with mode $i$ from $z_{k+1}$.


Under the lifted dynamics \eqref{eq::open_dynamics}, we formulate the following robust optimal control problem. Let $\mathbf{x}, \mathbf{z}, \mathbf{u}$ denote the history of the state $x_k$, lifted state $z_k$ and control $u_k$ over a finite horizon $\{0, \ldots, N\}$. Define the cost functional
\begin{equation}
    \begin{split}
    J\left(\mathbf{x},\mathbf{z},\mathbf{u}\right) &= \norm{x_{N}}_{P_{N}}^{2} 
        \\
        &+\sum_{k=0}^{N-1} \norm{x_{k}}^{2}_{Q_{\hat{\theta}_{k}}} + \norm{u_{k}}^{2}_{R_{\hat{\theta}_{k}}} + \norm{z_{k+1}}^{2}_{S_{\hat{\theta}_{k}}},
    \end{split}
\end{equation}
with $Q_{\hat{\theta}_{k}}, R_{\hat{\theta}_{k}}, S_{\hat{\theta}_{k}}\succ0$ and $P_{N}\succeq0$. The lifted state $z_{k+1}$ is weighted to account for all possible future state evolutions. The problem we aim to solve is then given by
\begin{equation} \label{prob::optmization_full}
    \min_{\mathbf{x},\mathbf{z},\mathbf{u}} \max_{\mathbf{C}^{\delta},\mathbf{q}^{\delta}} \mathbb{E}\left\{ J\left(\mathbf{x},\mathbf{z},\mathbf{u}\right) \Big| \hat{\theta}_{0},\hat{\theta}_{1},..., \hat{\theta}_{N}\right\}, ~~\text{sub. to dynamics } \eqref{eq::open_dynamics},
\end{equation}
where $\mathbf{C}^{\delta} = \{C_{\hat{{\theta}}_k}^{\delta}\}_{k=1}^{N-1}$ and $\mathbf{q}^{\delta} = \{ q_{\hat{{\theta}}_k\hat{\theta}_{k+1}}^{\delta}\}_{k=1}^{N-1}$. 

We now obtain a recursive formulation of \eqref{prob::optmization_full} by considering the associated Bellman's equation \cite{bertsekas2011dynamic}, 
\begin{align}\label{eq::full_optimization_recursive}
    &\min_{x_{k+1},z_{k+1},u_{k}} \max_{ C_{\hat{\theta}_{k}}^{\delta}, q_{\hat{{\theta}}_k\hat{\theta}_{k+1}}^{\delta}} J_{k} 
    ,
    ~~ \text{sub. to dynamics } \eqref{eq::open_dynamics},
\end{align}
with 
\begin{equation*}
    J_{k}=\norm{x_{k}}^{2}_{Q_{\hat{\theta}_{k}}} + \norm{u_{k}}^{2}_{R_{\hat{\theta}_{k}}} + \norm{z_{k+1}}^{2}_{S_{\hat{\theta}_{k}}}+ \mathbb{E}\{ J_{k+1}(x_{k+1},\hat{\theta}_{k+1}) | \hat{\theta}_{k}\},
\end{equation*}
which formulates the problem as a one-step optimization to be solved in a backward recursion.

For simplicity, we now denote $l = \hat{\theta}_k$ and $d = \hat{\theta}_{k+1}$. As the uncertain quantities $C_{l}^\delta$ and $q_{ld}^\delta$ are both dependent on the unknown true mode $i \in L_l$, an exact worst case formulation would require a joint maximization over the finite set of admissible modes within each cluster. In Section \ref{sec::RecursiveRegulation}, these uncertainties are embedded into structured norm-bounded uncertainty models, yielding a tractable robust relaxation and a closed-form solution.

{In addition to deriving a solution to the proposed control problem, we are concerned with the convergence of the resulting solution to a cluster-dependent stationary solution, enabling the analysis of stability of the closed-loop system}. To formally characterize this stability property, we adapt the notion of mean-square stability (MSS) from \cite{costa2005discrete} as follows.

\begin{definition}
System \eqref{eq::open_dynamics} is said to be MSS if, for any initial condition $(x_{0},\hat{\theta}_{0})$ and admissible uncertainty, $\mathbb{E}\{ \norm{x_{k}}^{2} \} \rightarrow 0$ as $k\rightarrow \infty$.
\end{definition}

The next section presents some auxiliary results in regularized least mean-square optimization. These results will be used in the solution of \eqref{eq::full_optimization_recursive} presented in Section~\ref{sec::RecursiveRegulation}.
\section{Auxiliary Lemmas} \label{sec::preliminaries}
This section consists of auxiliary results that are used to solve the optimization problem \eqref{eq::full_optimization_recursive}. Consider the regularized least mean-square optimization problem subject to uncertainties in the data and in the regularization weighting matrix 
\begin{equation}\label{eq::minimax_01}
\min_{\mathcal{x}} \max_{\delta \mathcal{A},\delta \mathcal{S}} {f(\mathcal{x})},
\end{equation}
with $f(\mathcal{x})=\norm{\mathcal{x}}_{\mathcal{S}+ \delta \mathcal{S}}^{2} + \norm{\left( \mathcal{A} + \delta \mathcal{A} \right)\mathcal{x} - \mathcal{b} }_{\mathcal{W}}^{2}$ and where $\mathcal{A}$ is a known matrix; $\mathcal{b}$ and $\mathcal{x}$ are the measurement and the minimization variables.  $\mathcal{S} \succ 0$ and $\mathcal{W} \succ 0$ are weighting matrices. The variables $\delta \mathcal{A}, \delta \mathcal{S}$ are uncertainties modeled by 
\begin{equation}\label{eq::minimax_02}
        \delta \mathcal{A}   =
    \mathcal{H} \Delta_{\mathcal{A}}
        E_{\mathcal{A}}, \ \ \ \delta \mathcal{S} = E_{\mathcal{S}}^T\Delta_{\mathcal{S}} E_{\mathcal{S}}.
\end{equation}
where $E_{\mathcal{A}}$ and $E_{\mathcal{S}}$ are known, $\mathcal{H}$ is a known nonzero matrix and $\Delta_{\mathcal{A}}$, $\Delta_{\mathcal{S}}=\Delta_{\mathcal{S}}^T$ are contractions with appropriate dimensions. 

The solution of the problem modeled in \eqref{eq::minimax_01} and \eqref{eq::minimax_02} is presented in the next lemma, and a proof is presented in Appendix~\ref{app:A}, following the framework in\cite{sayed2001framework}. 
\begin{lemma}\label{lemma::minimax_solution}
    Consider the min-max optimization problem (\ref{eq::minimax_01})-(\ref{eq::minimax_02}). Define $\bar{\mathcal{S}}=\mathcal{S}+E_{\mathcal{S}}^{T}E_{\mathcal{S}}$ and $\lambda^{\star}=\arg \inf_{\lambda >\norm{\mathcal{H}^{T}\mathcal{W}\mathcal{H}}} \Gamma(\lambda)$, where
    \begin{equation*}
        \begin{split}
            \Gamma(\lambda) &= \mathcal{b}^T\mathcal{W}_{\lambda}\mathcal{b}  \\ &-\mathcal{b}^T \mathcal{W}_{\lambda} \mathcal{A} \left( \bar{\mathcal{S}} + \lambda E_{\mathcal{A}}^T E_{\mathcal{A}}   +\mathcal{A}^{T}\mathcal{W}_\lambda \mathcal{A}\right)^{-1}\mathcal{A}^{T}\mathcal{W}_{\lambda}\mathcal{b},
        \end{split}
    \end{equation*}
    and
    \begin{equation}\label{eq::W_lambda}
    \mathcal{W}_{\lambda} = \mathcal{W} + \mathcal{W}\mathcal{H} \left( \lambda I -\mathcal{H}^{T}\mathcal{W}\mathcal{H} \right)^{-1} \mathcal{H}^{T}\mathcal{W}.
    \end{equation}
    
    Then the unique optimal solution and corresponding optimal cost $(\mathcal{x}^*, f(\mathcal{x}^{\star}))$ are given by
    \begin{equation*} 
        \begin{split}
           &(i) \ \mathcal{x}^{\star} = \left( \bar{\mathcal{S}} + \lambda^{\star} E_{\mathcal{A}}^{T}E_{\mathcal{A}} +\mathcal{A}^{T}\mathcal{W}_{\lambda^{\star}} \mathcal{A} \right)^{-1} \mathcal{A}^{T}\mathcal{W}_\lambda^{\star} \mathcal{b} ,
           \\&(ii) \ f(\mathcal{x}^{\star})=\Gamma(\lambda^\star).
       \end{split}
    \end{equation*}
\end{lemma}

\normalsize

\begin{remark}\label{remark::lambda_approx}
   As discussed in the work of \cite{SAYED2002361, bueno2022regulation}, the solution of problem $\lambda^{\star} \in \arg \inf_{\lambda > \norm{ \mathcal{H}^T\mathcal{W}\mathcal{H}}}  \Gamma(\lambda) $ is typically close to its lower bound, where $\beta>1^+$ is a good practical approximation for defining $\lambda^{\star}\approx \beta \norm{ \mathcal{H}^T\mathcal{W}\mathcal{H}}$. 
\end{remark}

The following auxiliary result, also proved in Appendix~\ref{app:A}, provides an upper bound for the maximization of a quadratic function subject to norm-bounded uncertainty. This result will be useful in developing the stability analysis in Section~\ref{sec::Stability}.
\begin{lemma}\label{lemma::maximization}
Consider $\mathcal{C}^{\delta}=\mathcal{C}+\delta\mathcal{C}$, with $\delta\mathcal{C}=\mathcal{H}\Delta E_{\mathcal{C}}$. Let $\mathcal{H}$ and $E_{\mathcal{C}}$ be known matrices, $\Delta$ a contraction and $\mathcal{W}\succ0$. Then, the following inequality holds
\begin{equation}\label{eq::maximization_bound}
 (\mathcal{C}+\delta\mathcal{C})^T\mathcal{W}(\mathcal{C}+\delta\mathcal{C}) \preceq \mathcal{C}^T\mathcal{W}_{\lambda}\mathcal{C} + \lambda E_{\mathcal{C}}^{T}E_{\mathcal{C}},
\end{equation}
for any $\lambda >\norm{\mathcal{H}^{T}\mathcal{W}\mathcal{H}}$, where $\mathcal{W}_{\lambda}$ is defined in \eqref{eq::W_lambda}.
\end{lemma}
Using these results, we are now ready to return to the setup formulated in Section~\ref{sec::problem_stat}.
\section{Robust Recursive Regulation}\label{sec::RecursiveRegulation}

We now discuss how a tractable solution is obtained to the recursive optimization problem \eqref{eq::full_optimization_recursive}. We first parameterize the uncertainties $q_{ld}^{\delta}$ and $C_l^{\delta}$ through structured bounded models
\begin{equation}
q_{ld}^{\delta}
= \bar{q}_{ld} + \alpha_{ld} \Delta_{ld}; ~~
C_l^{\delta}
= \bar{C}_l + M_{C_l}\Delta_{C_l}E_{C_l},
\label{eq::structured_uncertainties}
\end{equation}
where $\bar{q}_{ld}$ and $\bar{C}_l$ denote their respective nominal values, while $\alpha_{ld}\geq0$, non-zero $M_{C_l}$, and $E_{C_l}$ are known parameters defining the uncertainty structures, with $\Delta_{ld}$ and $\Delta_{C_l}$ denoting contractions of appropriate dimensions. This uncertainty description relaxes the original constraints on $q^{\delta}_{ld}$ and $C^{\delta}_{l}$ which required these quantities to lie on specific vertices.

Evaluating the conditional expectation within \eqref{eq::full_optimization_recursive},
\begin{equation}
    \mathbb{E}\{ J_{k+1}(x_{k+1},\hat{\theta}_{k+1}) | \hat{\theta}_{k}=l \} = x_{k+1}^T\Psi_{l,k+1}^{\delta}x_{k+1},
\end{equation}
where
\begin{equation}
    \Psi_{l,k+1}^{\delta} = \mathbb{E}\{ P_{d,k+1} | \hat{\theta}_{k}=l \} = \sum_{d \in \hat{\Theta}} \left( \bar{q}_{ld} + \alpha_{ld} \Delta_{ld}\right) P_{d,k+1}.
\end{equation}
The matrix $\Psi_{l,k+1}^{\delta}$ can then be decomposed as
\begin{equation}\label{eq:psi_factorization}
    \Psi_{l,k+1}^{\delta} = \bar{\Psi}_{l,k+1} + \delta \Psi_{l,k+1},
\end{equation}
with $\bar{\Psi}_{l,k+1}=\sum_{d \in \hat{\Theta}} \bar{q}_{ld} P_{d,k+1}$. For convenience, we factorize $\delta \Psi_{l,k+1}$ as follows:
\begin{equation}\label{eq::psi_factorization_02}
    \begin{split}
    &\delta \Psi_{l,k+1} = \mathcal{P}_{k+1}^{1/2^T} E_{Q_l}^{1/2^T}\Delta_Q E_{Q_l}^{1/2} \mathcal{P}_{k+1}^{1/2},
    \\
        &\mathcal{P}_{k+1}^{1/2^T}=\begin{bmatrix} P_{1,k+1}^{1/2^T} & ... & P_{h,k+1}^{1/2^T}\end{bmatrix},
        \\
        &E_{Q_l}^{1/2} = \text{blkdiag}(\alpha_{l1}^{1/2}I_{n},\dots,\alpha_{lh}^{1/2}I_{n}),
        \\
        &\Delta_Q = \text{blkdiag}(\Delta_{l1}I_{n},\dots,\Delta_{lh}I_{n}).
    \end{split}
\end{equation}
This factorization is well defined since $P_{d,k+1}\succeq0$ and $\alpha_{ld}\geq0$, ensuring the existence of the corresponding square roots.

We incorporate the dynamic constraints \eqref{eq::open_dynamics} into the recursive optimization problem through quadratic penalty terms \cite{cerri2017recursive,bueno2022regulation}, which will facilitate the derivation of a recursive solution. Let $\mu_1,\mu_2>0$ be penalty parameters and define
\begin{align}
    g(x_{k+1},z_{k+1},u_{k},x_{k}) &= \norm{z_{k+1} - Fx_{k} - Gu_{k}}_{\mu_1 I_{ns}}^{2}  \nonumber \\
    &+ \norm{C^{\delta}_{l}  z_{k+1} - x_{k+1}}_{\mu_2 I_{n}}^{2}.
\end{align}
\normalsize
Problem \eqref{eq::full_optimization_recursive} can then be rewritten as the unconstrained one-step min-max problem
\begin{equation}\label{eq::unconstrained_min_recursive}
    \begin{split}
        \min_{x_{k+1},z_{k+1},u_{k}} \max_{C_{l}^{\delta},q_{ld}^{\delta}} \left\{ \norm{x_{k}}^{2}_{Q_{l}} + \norm{u_{k}}^{2}_{R_{l}} + \norm{z_{k+1}}^{2}_{S_{l}}\right.\\\left.+\norm{x_{k+1}}_{\Psi_{l,k+1}^{\delta}}^{2} + g(x_{k+1},z_{k+1},u_{k},x_{k})
            \right\}.
    \end{split}
\end{equation}

We now cast \eqref{eq::unconstrained_min_recursive} into the framework of Lemma~\ref{lemma::minimax_solution} through the following identifications:
\begin{equation}\label{eq::identifications_beginning}
    \begin{gathered}
        \mathcal{A} \leftarrow 
        \begin{bmatrix}
        0 & 0 & 0
        \\
        0 & I_{ns} & -G
        \\
        -I_{n} & \bar{C}_{l} & 0
        \end{bmatrix},
        \ \ \ 
        \delta \mathcal{A} \leftarrow 
        \begin{bmatrix}
        0 & 0 & 0
        \\
        0 & 0 & 0
        \\
        0 & \delta C_{l} & 0
        \end{bmatrix},
        \\
        \mathcal{x}\leftarrow 
        \begin{bmatrix}
        x_{k+1}
        \\
        z_{k+1}
        \\
        u_{k}
        \end{bmatrix},
        \ \
        \mathcal{b} \leftarrow 
        \begin{bmatrix}
        I_{n}
        \\
        F
        \\
        0
        \end{bmatrix}x_{k},
        \\
        \mathcal{S} \leftarrow 
        \begin{bmatrix}
        \bar{\Psi}_{l,k+1} & 0 & 0 
        \\
        0 & S_{l} & 0
        \\
        0 & 0 & R_l
        \end{bmatrix},
        \ \ 
        \mathcal{W} \leftarrow 
        \begin{bmatrix}
        Q_l & 0 & 0
        \\
        0 & \mu_1 I_{ns} & 0
        \\
        0 & 0 & \mu_2 I_{n}
        \end{bmatrix},
    \end{gathered}
\end{equation}
with the corresponding uncertainty structures
\begin{equation}\label{eq::identifications_endings}
    \begin{gathered}
    \mathcal{H} \leftarrow 
    \begin{bmatrix}
    0
    \\
    0
    \\
    M_{C_l}
    \end{bmatrix}, 
    \ \ 
    \Delta_{\mathcal{A}} \leftarrow 
    \begin{bmatrix}
    \Delta_{C_l}
    \end{bmatrix},
    \ \
    E_{\mathcal{A}} \leftarrow 
    \begin{bmatrix}
    0 & E_{C_l} & 0
    \end{bmatrix},
    \\
    E_{\mathcal{S}} \leftarrow 
    \begin{bmatrix}
    E_{Q_l}^{1/2}\mathcal{P}_{k+1}^{1/2} & 0 & 0
    \end{bmatrix},
    \ \
    \Delta_{\mathcal{S}} \leftarrow 
    \begin{bmatrix}
        \Delta_{Q_l}
    \end{bmatrix}.
    \end{gathered}
\end{equation}
Under these identifications, \eqref{eq::unconstrained_min_recursive} is identical to \eqref{eq::minimax_01}, and we can thus use Lemma~\ref{lemma::minimax_solution} to obtain the minimizing variables in closed-form. This closed-form solution results in a cluster-dependent feedback gain when we enforce the system dynamics exactly by letting $\mu_1 \to \infty$ but keeping $\mu_2 > 0$ fixed. We note that this distinction reflects the differing roles of $\mu_1$ and $\mu_2$: $\mu_1$ is associated with a relation that must hold exactly whereas $\mu_2$ penalizes a relation affected by mode uncertainty and  thus remains finite.

The finite-horizon solution is obtained through the backward recursion indexed through $\{N-1,\ldots,0\}$. However, these indices are used only to encode the ordering of successive iterates in implementing the resulting recursion. For a convenient computational representation, we reverse the order of indexing to $\{0, \ldots, N-1\}$. The cluster-dependent feedback gain is then derived using the following theorem, which is proved in Appendix~\ref{app:B}.

\begin{theorem}\label{thm::robust_riccati}
Consider the one-step optimization problem \eqref{eq::unconstrained_min_recursive} associated with the observed cluster $l=\hat{\theta}_k$ at time $k$. Let $\mu_2>0$ and $\beta>1$. For a set of weighting matrices $S_{l},  Q_{l},  R_{l}, \succ0$ and a collection of value matrices $P_{d,k}\succeq0$, $d\in\hat{\Theta}$, define the matrix
\begin{equation}\label{eq::Psi_theorem}
    \Psi_{l,k}
    =
    \sum_{d\in\hat{\Theta}}
    \left(
        \bar{q}_{ld}+\alpha_{ld}
    \right)P_{d,k},
\end{equation}
and auxiliary functions
\begin{equation}\label{eq::aux_functions}
    \begin{split}
        &\lambda_l
        =
        \beta\mu_2
        \left\|M_{C_l}^{T}M_{C_l}\right\|, W_{l}^{C}=\left(\mu_2^{-1}I_n-\lambda_l^{-1}M_{C_l}M_{C_l}^{T}\right)^{-1},
        \\
        &\Omega_{l,k}=S_l+\lambda_l E_{C_l}^{T}E_{C_l}
        \\
        & \hspace{1cm} +\bar{C}_l^{T}\left(W_{l}^{C}-W_{l}^{C}\left(\Psi_{l,k}+W_{l}^{C}\right)^{-1}W_{l}^{C}\right)\bar{C}_l.  
    \end{split}
\end{equation}
The optimal cost and the minimizing control action are $J_k^{\star}(x_k,l)=x_k^{T}P_{l,k}x_k$ and $u_k^{\star}=K_{l,k}x_k$, where $P_{l,k}$ is recursively given by
\begin{equation}\label{eq::robust_riccati}
\begin{split}
    &P_{l,k+1}=Q_l+F^{T}\Omega_{l,k}F 
    \\
    &-F^{T}\Omega_{l,k}G\left(R_l+G^{T}\Omega_{l,k}G\right)^{-1}G^{T}\Omega_{l,k}F ,
\end{split}
\end{equation}
and the cluster dependent feedback gain as
\begin{equation}\label{eq::robust_gain}
    K_{l,k}=-\left(R_l+G^{T}\Omega_{l,k}G\right)^{-1}G^{T}\Omega_{l,k}F,
\end{equation}
\end{theorem}


\begin{remark}
Although the uncertain transition probabilities $q_{ld}^{\delta}$ in \eqref{eq::structured_uncertainties} must satisfy condition \eqref{eq::prob_Markov_cluster_02}, the coefficients $\bar q_{ld}+\alpha_{ld}$ in $\Psi_{l,k}$ \eqref{eq::Psi_theorem} are not bound by this restriction. Given that $\sum_{d\in\hat{\Theta}}(\bar q_{ld}+\alpha_{ld}) = 1+\sum_{d\in\hat{\Theta}}\alpha_{ld}$, these coefficients emerge from the robust maximization of transition uncertainty and must not be treated as valid transition probabilities. Furthermore, $\Psi_{l,k+1}$  functions strictly as a robust upper bound weighting within the recursion, rather than as the conditional expectation of a modified Markov chain. 
\end{remark}
\section{Convergence and Stability}\label{sec::Stability}
We now prove that the sequence of matrices $\left\{\left\{ P_{l,k} \right\}_{l \in \hat{\Theta}} \right\}_{k=0}^{\infty}$ generated by \eqref{eq::robust_riccati} converges to a stationary solution $P_l$ for each cluster $l\in\hat{\Theta}$ as $k\to\infty$. Using this property, we analyze closed-loop stability under the stationary cluster-dependent control law that is subsequently obtained \eqref{eq::robust_gain}.




\subsection{Convergence}
\label{subsec::convergence}
The recursion \eqref{eq::robust_riccati} generates a sequence $\{P_{l,k}\}_{k \geq 0}$ for each cluster $l\in\hat{\Theta}$. The following proposition, proved in Appendix~\ref{app:C}, establishes that these sequences remain bounded and converge to finite limiting matrices $P_l$.


\begin{proposition}\label{prop::convergence}
    Assume fixed $S_l,Q_l, R_l \succ 0$, $\mu_{2}>0$, and $\lambda_l>\mu_{2}\norm{M_{C_l}^TM_{C_l}}$ such that $W^{C}_{l}\succ 0$ for all $l \in \hat{\Theta}$. Initialize the coupled Riccati recursion \eqref{eq::robust_riccati} with initial condition $P_{l,0}=0$ for all $l \in \hat{\Theta}$. Then the sequence $\left\{\left\{ P_{l,k} \right\}_{l \in \hat{\Theta}} \right\}_{k=0}^{\infty}$ is
    \begin{enumerate}
        \item Monotonically non-decreasing,
        \item Bounded above, and
        \item Convergent to a finite set $\{ P_{l} \}_{l\in \hat{\Theta}}$.
    \end{enumerate}
\end{proposition}

An important feature of Proposition \ref{prop::convergence} is that no additional stabilizability assumption is required to establish convergence of the Riccati sequence. We are now in a position to analyze stability of the closed-loop system.



\subsection{Stability}\label{subsec::stability}

To establish the stability of system \eqref{eq::open_dynamics} under the uncertainties in \eqref{eq::structured_uncertainties} and the cluster-dependent control law $\{K_l\}_{l\in\hat{\Theta}}$, we consider the coupled quadratic Lyapunov function from \cite{costa2005discrete},
\begin{equation*}
    V_{k}\left( x_{k},l=\hat{\theta}_{k} \right)=x_{k}^T P_{l}x_{k}, \ \ \ P_{l} \succ 0.
\end{equation*}
Since the Lyapunov difference depends on admissible uncertainties, we first derive an upper bound on the difference $V_{k+1}-V_k$. The converged Riccati equation is then used to establish conditions under which this upper bound is negative definite.

Under the closed-loop feedback $u_k=K_l x_k$, the dynamics in \eqref{eq::open_dynamics} is given by $x_{k+1}=C_{l}^\delta \left( F + G K_{l} \right) x_{k}$ and 
\begin{align}
    V_{k+1}-V_{k} &= \mathbb{E}\{x_{k+1}^T P_{d}x_{k+1}\}-  x_{k}^T P_{l} x_{k} \nonumber \\
    &= x_{k}^T \left( F_{CL}^T C_{l}^{\delta^T} \Psi_{l}^{\delta}  C_{l}^\delta F_{CL} -  P_{l} \right) x_{k},\label{eq::Lyapunov_Difference}
\end{align}
where $F_{CL}=F+GK_l$, and $\Psi_l^\delta$ denotes the conditional expected-value matrix defined in \eqref{eq:psi_factorization}. The matrix in \eqref{eq::Lyapunov_Difference} is bounded in the following lemma, proved in Appendix~\ref{app:C}.


\begin{lemma}\label{lemma::upper_bound_Lyapu}
    For $\lambda_l > \norm{M_{C_l}^T \Psi_{l}M_{C_l}}$, it follows that 
    \begin{equation}\label{eq::upper_bound_Lyapunov}
        \begin{split}
            & F_{CL}^T C_{l}^{\delta^T} \Psi_{l}^{\delta}  C_{l}^\delta F_{CL} -  P_{l}
            \\
            \preceq~ &F_{CL}^T \bar{C}_{l}^T \Psi_{l}^{C} \bar{C}_{l}F_{CL} + \lambda_l F_{CL}^TE_{C_l}^T E_{C_l}F_{CL} - P_l,
        \end{split}
    \end{equation}
    with $\Psi_{l}^{C} = \Psi_l + \Psi_l M_{C_l} \left( \lambda_l I - M_{C_l}^T \Psi_{l}M_{C_l}\right)^{-1} M_{C_l}^T \Psi_l$.
\end{lemma}

Using this bound, we establish mean-square stability of the closed-loop system in the following theorem, which is proved in Appendix~\ref{app:C}.


\begin{theorem}\label{thm::stability}
    Assume fixed $S_l,Q_l, R_l \succ 0$, $\mu_{2}>0$, and $\lambda_l>\mu_{2}\norm{M_{C_l}^TM_{C_l}}$ such that $W^{C}_{l}\succ 0$ for all $l \in \hat{\Theta}$. If
    \begin{equation}
        \begin{split}
        &\lambda_l>\norm{M_{C_l}^T \Psi_{l}M_{C_l}},
        \\
        &S_{l} \succeq \bar{C}_l^{T}\left(\Psi_{l}^{C} - W^{C}_l+W^{C}_l\left(\Psi_{l}+W^{C}_l\right)^{-1}W^{C}_l\right) \bar{C}_{l},
        \end{split}
    \end{equation}
    for all $l \in \hat{\Theta}$, then the close-loop system dynamics $x_{k+1}=C_{l}^{\delta}F_{CL}x_{k}$ is MSS under feedback $\{ K_{l} \}_{l\in \hat{\Theta}}$ derived from the Riccati solution in Theorem \ref{thm::robust_riccati} when subjected to uncertainty modeled as in \eqref{eq::structured_uncertainties}.
\end{theorem}

This section concludes the analysis of our solution proposed in Section~\ref{sec::RecursiveRegulation}, to the problem posed in Section~\ref{sec::problem_stat}. {We remark that the result presented in Theorem \ref{thm::stability} guarantees stability, though it may have a certain level of conservatism as conditions \eqref{eq::structured_uncertainties} may not always be satisfied}. In the next section, we illustrate the performance of our solution in the following section on the model of an unmanned aerial vehicle.
\section{Example Revisited} \label{sec::Examples}
\label{subsec::UAV_example}

Consider the fault-susceptible lateral-directional dynamics of an unmanned aerial vehicle presented in \cite{todorov2018detector}, with the nominal dynamics $A_1,B_1$ from \cite[Example 6.1]{todorov2018detector} denoted by $A_{OK},B_{OK}$. The mode transitions are given in Example 1, and the uncertain transition probabilities between clusters $A$ and $B$ are given by
\begin{equation}
    \mathbb{L}^{\delta} =
    \begin{bmatrix}
        0.95 + 0.05 \Delta & 0.05 + 0.05(- \Delta)
        \\
        1.00 & 0.00
    \end{bmatrix}.
\end{equation}
Thus, $\bar q_{AA}=0.95$, $\bar q_{AB}=0.05$, $\bar q_{BA}=1$, and $\bar q_{BB}=0$, with uncertainty bounds $\alpha_{AA}=\alpha_{AB}=0.05$ and $\alpha_{BA}=\alpha_{BB}=0$. This representation preserves the stochasticity of the transition matrix. The output matrices selecting the state from the lifted state $z_{k+1}$ are modeled as in \eqref{eq::structured_uncertainties}, and are given by
\begin{equation}
    \bar{C}_{A} = 
    \begin{bmatrix}
        \frac{1}{2}I_{4} & \frac{1}{2}I_{4} & 0
    \end{bmatrix},
    \ \ \ 
    \bar{C}_{B} = 
    \begin{bmatrix}
        0 & 0 & I_{4}
    \end{bmatrix},
\end{equation}
\begin{equation}
    E_{C_A} = 
    \begin{bmatrix}
        \frac{1}{2}I_{4} & -\frac{1}{2}I_{4} & 0    
    \end{bmatrix},
    \ \ \
    E_{C_B} = 
    \begin{bmatrix}
        0 & 0 & 0    
    \end{bmatrix},
\end{equation}
with $M_{C_A}=M_{C_B}=I_{4}$.

{The fault modes represent partial loss of actuator effectiveness, with $B_{UF}=B_{DF}=0.1B_{OK}$. Alongside this model, we introduce an additional degradation factor $s\geq 1$ in the state matrix, defining $A_{UF}=A_{DF}=sA_{OK}$. This factor is used for numerical simulation comparisons; by increasing it, we consider faulty systems whose dynamics deviate from their nominal values.} 

We implement the procedure in Theorem~\ref{thm::robust_riccati}\footnote{https://github.com/leadcatlab/LumpableClusterControl} by setting $P_{0}=I_{4}$, $Q_{A}=Q_{B}=I_{4}$, $R_{A}=R_{B}=I_{2}$, $S_{A}=S_{B}=I_{12}$, with $\mu_2=100$, and $\beta=1.0001$. {We employ results from \cite{do2002h2} and \cite{todorov2018detector} as comparative baselines}. All controllers are designed using identity state and control weighting matrices. The closed-loop comparison is performed without external disturbances, using the same initial condition and the same realization of the Markov mode sequence. Furthermore, the comparison methods were implemented in Python using the \texttt{CVXPY} optimization framework.

\begin{figure}[h!]
    \centering
    \includegraphics[width=0.95\linewidth]{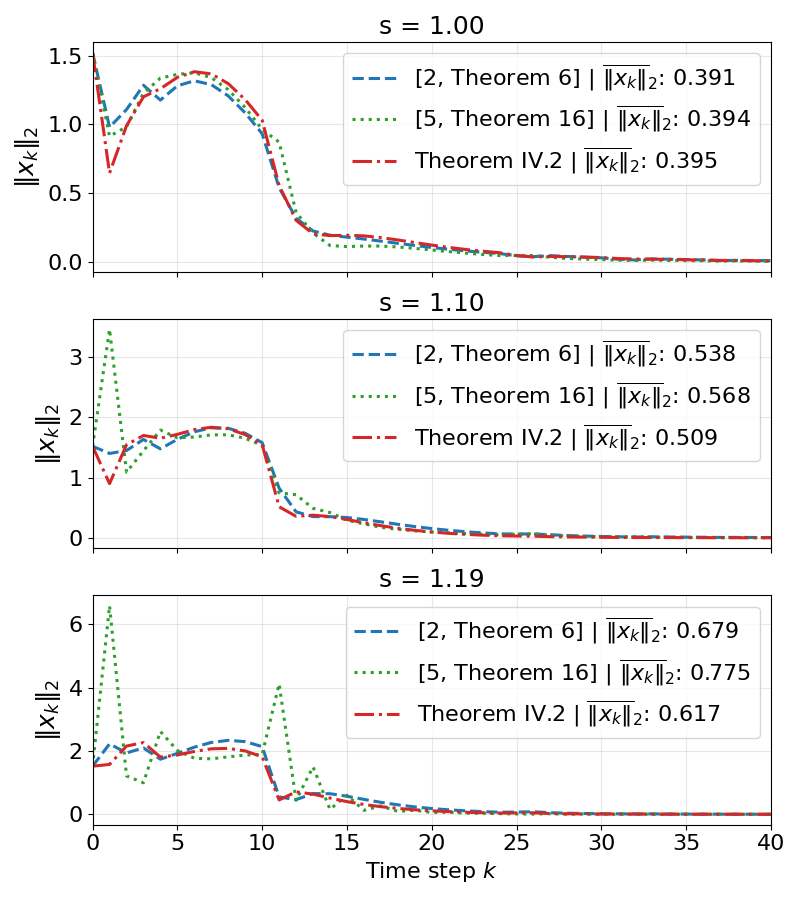}
    \caption{Closed-loop state norm $\norm{x_k}_2$ for the proposed controller in Theorem \ref{thm::robust_riccati} and comparison methods, under the same mode realization and initial condition. The parameter $s$ represents increasing degradation of the system matrix $A_{OK}$ in the fault modes. The value reported in the legend correspond to the average state norm over the simulation time.}
    \label{fig::Sim}
\end{figure}

In Figure \ref{fig::Sim} we show results for three different degradation factor $s=\{ 1.00, 1.10,1.19 \}$. For $s=1.00$, the three controllers exhibit similar closed-loop behavior, whereas their differences become more pronounced as the degradation factor increases. In particular, the proposed controller presents a less severe transient response as $s$ increases. 

The interpretation of these results is twofold. First, the min-max formulation explicitly accounts for worst-case mode transitions and dynamics, introducing a degree of conservatism that becomes beneficial as the discrepancy among the modes within a cluster increases. Second, the proposed solution ultimately reduces to a closed-form Riccati recursion, whose convergence is established in Proposition \ref{prop::convergence}. Thus, its implementation consists of directly evaluating the Riccati recursion without relying on general-purpose optimization solvers, considerably reducing computational effort.


\section{Conclusion}\label{sec::Conclusion}

In this paper we address the {problem of} controlling Markov jump linear systems with lumpable cluster observations. We {pose this problem as} a min-max quadratic optimization problem that maximizes the worst-case mode associated uncertainty while choosing a controller to minimize the resulting cost. We provide a lemma to address this problem by modifying a well-known result from the literature. Under further considerations on the uncertainty model, we derive a closed-form solution based on a Riccati recursion. We demonstrate the numerical tractability of the solution through a convergence analysis, and provide stability guarantees to the convergent solution. Finally, we provide a numerical example from the literature, with benchmark comparisons. {Future work will focus on constructing the uncertain cluster transition matrix from the original Markov chain and relaxing lumpability requirements.}

\appendices
\section{} \label{app:A}
\subsection{Proof of Lemma \ref{lemma::minimax_solution}}
First, the maximization problem over $\delta \mathcal{S}$ is carried out. It can be stated as,
\begin{equation*}
    \max_{\delta \mathcal{A}} \left\{ \norm{r}_{\mathcal{W}}^{2}  + \norm{\mathcal{x}}_{\mathcal{S}}^{2} + \max_{\Delta_\mathcal{S}} \norm{E_{\mathcal{S}}\mathcal{x}}_{\Delta_{\mathcal{S}} }^{2}\right\},
\end{equation*}
where $r=\left( \mathcal{A} + \delta \mathcal{A} \right)\mathcal{x} - \mathcal{b} $. By employing the Rayleigh inequality (F3), we have that for all $\mathcal{x}$, $\rho_{min}(\Delta_{\mathcal{S}})\norm{E_{\mathcal{S}}\mathcal{x}}^2 \leq \norm{E_{\mathcal{S}}\mathcal{x}}^2_{\Delta_{\mathcal{S}}} \leq \rho_{max}(\Delta_{\mathcal{S}})\norm{E_{\mathcal{S}}\mathcal{x}}^2.$
Thus, since $\rho_{max}(\Delta_{\mathcal{S}})\leq1$, it follows that $\max_{\Delta_\mathcal{S}} \norm{E_{\mathcal{S}}\mathcal{x}}_{\Delta_{\mathcal{S}} }^{2}=\norm{\mathcal{x}}_{E_{\mathcal{S}}^TE_{\mathcal{S}}}^2$. Using $\bar{\mathcal{S}}=\mathcal{S}+E_{\mathcal{S}}^{T}E_{\mathcal{S}}$, the problem reduces to the formulation considered in \cite{sayed2001framework ,cerri2017recursive}:
\begin{equation*}
    \min_{\mathcal{x}} \max_{\delta \mathcal{A}} \left\{ \norm{\mathcal{x}}_{\bar{\mathcal{S}}}^{2} + \norm{\left( \mathcal{A} + \delta \mathcal{A} \right)\mathcal{x} - \mathcal{b}}_{\mathcal{W}}^{2}  \right\},
\end{equation*}
whose solution $\mathcal{x}^\star$ follows immediately from \cite[Theorem 1]{sayed2001framework}. Following the same reference, the optimum cost is
\begin{align*}
        &f(\mathcal{x}^{\star}) = \\ 
            &\hspace{0.2em} \mathcal{x}^{\star^T} \left(
            \bar{\mathcal{S}} + \lambda^{\star} E_{\mathcal{A}}^T E_{\mathcal{A}} \right)\mathcal{x}^{\star} 
            +(\mathcal{A}\mathcal{x}^{\star} - \mathcal{b})^T \mathcal{W}_{\lambda^{\star}} (\mathcal{A}\mathcal{x}^{\star} - \mathcal{b}),
\end{align*}
\normalsize
which on simplification yields the result in the lemma.

\subsection{Proof of Lemma \ref{lemma::maximization}}
Start with the contraction $\Delta^T\Delta -I\preceq 0$ and multiply left and right sides by $E_{\mathcal{C}}^T$ and $E_{\mathcal{C}}$ so that $E_{\mathcal{C}}^T \left( \Delta^T\Delta -I \right) E_{\mathcal{C}} \preceq 0$. Let $\mathcal{E}=\Delta E_{\mathcal{C}}$, such that $E_{\mathcal{C}}^T \left( \Delta^T\Delta -I \right) E_{\mathcal{C}} = \mathcal{E}^T\mathcal{E} -E_{\mathcal{C}}^TE_{\mathcal{C}}$ and $\delta \mathcal{C} = \mathcal{H}\mathcal{E}$, then
\begin{equation*}
    \begin{split}
        &(\mathcal{C}+\mathcal{H}\mathcal{E})^T\mathcal{W}(\mathcal{C}+\mathcal{H}\mathcal{E}) 
        \\
        &\preceq (\mathcal{C}+\mathcal{H}\mathcal{E})^T\mathcal{W}(\mathcal{C}+\mathcal{H}\mathcal{E}) -\lambda\left( \mathcal{E}^T\mathcal{E} -E_{\mathcal{C}}^TE_{\mathcal{C}}\right)=L,
    \end{split}
\end{equation*}
for any $\lambda > 0$. Rearranging,
$L = \mathcal{C}^T\mathcal{W}\mathcal{C} + \mathcal{C}^T\mathcal{W}\mathcal{H}\mathcal{E} + \mathcal{E}^T \mathcal{H}^T \mathcal{W}\mathcal{C} - \mathcal{E}^T \left( \lambda I - \mathcal{H}^T \mathcal{W} \mathcal{H}\right) \mathcal{E} + \lambda E_{\mathcal{C}}^T E_{\mathcal{C}}$. 
Define $M_{\lambda}=\lambda I - \mathcal{H}^T \mathcal{W} \mathcal{H}$. For $\lambda > \norm{\mathcal{H}^T \mathcal{W} \mathcal{H}}$ the matrix $M_{\lambda}\succ 0$. Completing squares with respect to $\mathcal{E}$ and $M_{\lambda}^{-1}\mathcal{H}^T\mathcal{W}\mathcal{E}$,
\begin{equation*}
    \begin{split}
        L = \mathcal{C}^T\mathcal{W}\mathcal{C} + \lambda E_{\mathcal{C}}^T E_{\mathcal{C}} + \mathcal{C}^T\mathcal{W}\mathcal{H}M_{\lambda}^{-1}\mathcal{H}^T\mathcal{W}\mathcal{C} \\- \left( \mathcal{E}-  M_{\lambda}^{-1}\mathcal{H}^T\mathcal{W}\mathcal{E}\right)^TM_{\lambda}\left( \mathcal{E} -  M_{\lambda}^{-1}\mathcal{H}^T\mathcal{W}\mathcal{E}\right).   
    \end{split}
\end{equation*}
As $M_{\lambda}\succ 0$, $L$ has an upper bound 
\begin{equation*}
    L \preceq \mathcal{C}^T\mathcal{W}\mathcal{C} + \lambda E_{\mathcal{C}}^T E_{\mathcal{C}} + \mathcal{C}^T\mathcal{W}\mathcal{H}M_{\lambda}^{-1}\mathcal{H}^T\mathcal{W}\mathcal{C}.
\end{equation*}
By defining $\mathcal{W}_{\lambda}$ as in \eqref{eq::W_lambda}, the result follows.
\section{} \label{app:B}
\subsection{Proof of Theorem \ref{thm::robust_riccati}}
Considering the identifications \eqref{eq::identifications_beginning}-\eqref{eq::identifications_endings} and Lemma \ref{lemma::minimax_solution}, the modified weighting matrix associated with the uncertainty in the transition probabilities is given by
$$
\bar{\mathcal S} = \mathcal S+E_{\mathcal S}^{T}E_{\mathcal S}
= \text{blkdiag}\left( \bar\Psi_{l,k} + \mathcal P_{k}^{1/2^T} E_{Q_l} \mathcal P_{k}^{1/2}, \,S_l,\,R_l \right).
$$
Since $\mathcal P_{k}^{1/2^T} E_{Q_l} \mathcal P_{k}^{1/2} = \sum_{d\in\hat{\Theta}}\alpha_{ld}P_{d,i},$ the first diagonal block of $\bar{\mathcal{S}}$ becomes $\Psi_{l,k} = \sum_{d\in\hat{\Theta}} \left(\bar q_{ld}+\alpha_{ld}\right)P_{d,k},$ as defined in (\ref{eq::Psi_theorem}). Furthermore, the scalar $\lambda_l$ introduced by Lemma~\ref{lemma::minimax_solution} must satisfy its corresponding lower-bound condition. As $\mathcal H^T\mathcal W\mathcal H=\mu_2M_{C_l}^{T}M_{C_l}$, following Remark~\ref{remark::lambda_approx}, we choose $\lambda_l=\beta\mu_2\left\|M_{C_l}^{T}M_{C_l}\right\|$ with $\beta>1^{+}$.

From Lemma   \ref{lemma::minimax_solution}(ii), the optimal cost function can, after some algebraic manipulation, be written as
\begin{equation*}
\begin{split}
P_{l,k+1}=Q_l+\mu_1 F^T F-
\begin{bmatrix}
0
\\
\mu_1 F
\\
-\mu_1 G^T F
\end{bmatrix}^{T}
\mathcal{M}_{l,k}
\begin{bmatrix}
0
\\
\mu_1F
\\
-\mu_1G^TF
\end{bmatrix},
\end{split}
\end{equation*}
where
\begin{equation*}
\mathcal{M}_{l,k}
=
\begin{bmatrix}
\Psi_{l,k}+W_{l}^{C} & -W_{l}^{C}\bar C_l & 0
\\
-\bar C_l^TW_{l}^{C} & \widetilde{\mathcal Z}_{l} & -\mu_1G
\\
0 & -\mu_1G^T & \widetilde R_{l}
\end{bmatrix}^{-1},
\end{equation*}
with $\widetilde{\mathcal Z}_{l}=S_l+\mu_1I_{ns}+\bar C_l^TW_{l}^{C}\bar C_l+\lambda_lE_{C_l}^TE_{C_l}$ and $\widetilde R_{l} = R_l+\mu_1G^TG$. Taking the Schur complement of
$\Psi_{l,k}+W_{l}^{C}$ gives
\begin{equation*}
P_{l,k+1} = Q_l + \mu_1 F^T F-
\begin{bmatrix}
\mu_1F
\\
-\mu_1 G^T F
\end{bmatrix}^{T}
\mathcal N_{l,k}
\begin{bmatrix}
\mu_1 F
\\
-\mu_1 G^T F
\end{bmatrix},
\end{equation*}
where
\begin{equation*}
\mathcal N_{l,k}
=
\begin{bmatrix}
\Omega_{l,k} + \mu_1 I_{ns} & -\mu_1G
\\
-\mu_1 G^T & \widetilde R_{l}
\end{bmatrix}^{-1},
\end{equation*}
and $\Omega_{l,k}$ is defined in (\ref{eq::aux_functions}). Introducing  the auxiliary matrices $\tilde{\Omega}_{l,k} = \Omega_{l,k} + \mu_1 I_{ns}$ and $\Omega_{l,k}^{\mu_{1}}=\mu_1 I_{ns}-\mu_1^2 \left(\mu_1 I_{ns}+ \Omega_{l,k}\right)^{-1}$, we use Property F2 in Section \ref{section:int} to rewrite $\Omega_{l,k}^{\mu_{1}}$ as $\Omega_{l,k}^{\mu_1} = \Omega_{l,k}\left( I_{ns} + \mu_{1}^{-1}\Omega_{l,k} \right)^{-1}$. Applying the Banachiewicz inversion \cite{tian2005schur} formula to $\mathcal N_{l,k}$ yields
\begin{equation*}
    \begin{split}
        \mathcal N_{11}=&\tilde{\Omega}_{l,k}^{-1}+\mu_1^2\tilde{\Omega}_{l,k}^{-1}G\left(R_l+G^T\Omega_{l,k}^{\mu_1}G\right)^{-1}G^T\tilde{\Omega}_{l,k}^{-1},
        \\
        \mathcal N_{12}=&\mu_1\tilde{\Omega}_{l,k}^{-1}G\left(R_l+G^T\Omega_{l,k}^{\mu_1}G\right)^{-1},
        \\\mathcal N_{21}\!=&\mu_1\left(R_l\!+\!G^T\Omega_{l,k}^{\mu_1}G\right)^{-1}G^T\tilde{\Omega}_{l,k}^{-1}, 
        \mathcal N_{22}\!=\!\left(R_l\!+\!G^T\Omega_{l,k}^{\mu_1}G\right)^{-1}.
    \end{split}
\end{equation*}
Substituting these blocks into the quadratic form above and collecting
terms yields the stated Riccati recursion.  The feedback gain follows directly from relation (i) of Lemma~\ref{lemma::minimax_solution}. In particular, 
\begin{equation*}
    \begin{bmatrix}
        x_{k+1}
        \\
        z_{k+1}
        \\
        u_{k}
    \end{bmatrix}
    =
    \mathcal{M}_{l,k}
    \begin{bmatrix}
    0
    \\
    \mu_1F
    \\
    -\mu_1G^TF
    \end{bmatrix}x_{k}; ~~~
    \begin{bmatrix}
        z_{k+1}
        \\
        u_{k}
    \end{bmatrix}
    =
    \mathcal{N}_{l,k}
    \begin{bmatrix}
    \mu_1F
    \\
    -\mu_1G^TF
    \end{bmatrix}x_{k},
\end{equation*}
which reduces  to
\begin{equation*}
    u_{k} = (\mu_{1}\mathcal{N}_{21}  - \mu_{1}\mathcal{N}_{22}G^T  )Fx_{k}.
\end{equation*}
Finally the optimal cost $J_k^{\star}(x_k,l)=x_k^{T}P_{l,k}x_k$ and control $u_k^{\star}=K_{l,k}x_k$, are given by
\begin{equation}
\begin{split}
    &P_{l,k+1}=Q_l+F^{T}\Omega_{l,k}^{\mu_{1}}F 
    \\
    &-F^{T}\Omega_{l,k}^{\mu_{1}}G\left(R_l+G^{T}\Omega_{l,k}^{\mu_{1}}G\right)^{-1}G^{T}\Omega_{l,k}^{\mu_{1}}F ,
    \\
    &K_{l,k}=-\left(R_l+G^{T}\Omega_{l,k}^{\mu_{1}}G\right)^{-1}G^{T}\Omega_{l,k}^{\mu_{1}}F.
\end{split}
\end{equation}
Taking the limit of $\Omega_{l,k}^{\mu_1} = \Omega_{l,k}\left( I_{ns} + \mu_{1}^{-1}\Omega_{l,k} \right)^{-1}$, $\lim_{\mu_{1} \rightarrow \infty}  \Omega_{l,k}^{\mu_1}  = \Omega_{l,k},$ thus proving the result.

\section{} \label{app:C}

\subsection{Proof of Proposition \ref{prop::convergence}}
The proof proceeds in three steps. (1) First, we establish monotonicity. (1.1) We show that initializing $P_{l,0}=0$ yields $P_{l,1}\succeq P_{l,0}$ for all $l \in \hat{\Theta}$. (1.2) For the induction step we show that $P_{d,k+1} \succeq P_{d,k}$ for all $d\in \hat{\Theta}$ leads to $P_{l,k+2} \succeq P_{l,k+1}$ for all $l \in \hat{\Theta}$. (2) Second, we show that there exists an upper bound to the sequence. (3) Concluding the proof, convergence follows from (1) and (2).

\textbf{(1.1)} Taking $P_{l,0}=0$, $\Psi_{l,0}=\sum_{d\in \hat{\Theta}} (q_{ld}+\alpha_{ld})P_{d,0}=0$. Hence $\Omega_{l,0} = S_{l} + \lambda_l E_{C_l}^T E_{C_l} \succ 0$, resulting in
\begin{equation*}
    P_{l,1} = Q_{l} + F^T\left( \Omega_{l,0}^{-1} + G R_{l}^{-1} G^T \right)^{-1}F \succ 0,
\end{equation*}
using Woodbury formula (F1). Therefore, for an initial condition $P_{l,0}=0$, $P_{l,1}\succeq P_{l,0}$.

\textbf{(1.2)} Assume that $P_{d,k+1} \succeq P_{d,k}$ for all $d\in \hat{\Theta}$., it follows that $\Psi_{l,k+1} - \Psi_{l,k} = \sum_{d \in \hat{\Theta}} (q_{ld}+\alpha_{ld})\left( P_{d,k+1} - P_{d,k} \right)$. Since $(q_{ld}+\alpha_{ld})\geq 0$, the inequality $\Psi_{l,k+1} \succeq \Psi_{l,k}$ holds and
\begin{equation*}
    \begin{split}
    &\Omega_{l,k+1}-\Omega_{l,k} = \bar{C}_{l}^T \left[ -W_{l}^{C}\left(\Psi_{l,k+1}+W_{l}^{C}\right)^{-1}W_{l}^{C} \right. \\
    &\left. + W_{l}^{C} \left(\Psi_{l,k}+W_{l}^{C}\right)^{-1}W_{l}^{C}\right] \bar{C}_{l}.
    \end{split}
\end{equation*}
Since $\Psi_{l,k+1} \succeq \Psi_{l,k}$, it follows that 
$\left(\Psi_{l,k+1} + W^{C}_{l}\right)^{-1} \preceq  \left( \Psi_{l,k}+ W^{C}_{l} \right)^{-1}$ and 
\begin{equation*}
    -W_{l}^{C}\left(\Psi_{l,k+1}+W_{l}^{C}\right)^{-1}W_{l}^{C} \succeq - W_{l}^{C}\left(\Psi_{l,k}+W_{l}^{C}\right)^{-1}W_{l}^{C},
\end{equation*}
which shows that $P_{d,k+1}\succeq P_{d,k}$  implies $\Omega_{l,k+1}\succeq \Omega_{l,k}$. Finally, using the Woodbury matrix identity,
\begin{equation*}
    \begin{split}
        &P_{l,k+2} - P_{l,k+1} = 
        \\
        &=F^T \left[ \left( \Omega_{l,k+1}^{-1} + GR_{l}^{-1}G^T \right)^{-1} - \left(\Omega_{l,k}^{-1} + GR_{l}^{-1}G^T \right)^{-1} \right] F.
    \end{split}
\end{equation*}
Consider that $\Omega_{l,k+1}\succeq \Omega_{l,k}$ implies $\Omega_{l,k+1}^{-1}\preceq \Omega_{l,k}^{-1}$, which results in $\Omega_{l,k+1}^{-1} + GR_{l}^{-1}G^T\preceq \Omega_{l,k}^{-1} + GR_{l}^{-1}G^T$, and therefore
\begin{equation*}
    \left(\Omega_{l,k+1}^{-1} + GR_{l}^{-1}G^T \right)^{-1} \succeq \left(\Omega_{l,k}^{-1} + GR_{l}^{-1}G^T\right)^{-1}.
\end{equation*}
Then, $P_{d,k+1} \succeq P_{d,k}$ implies $P_{l,k+2}\succeq P_{l,k+1}$ , for all $l,d \in \hat{\Theta}$.  

\textbf{(2)} Now, an upper bound on $P_{l,k+1}$ is defined. Since $W^{C}_{l}\succ0$ and $\Psi_{l,k}\succeq 0$ , 
\begin{equation*}
    \begin{split}
        &\Omega_{l,k}\!=\! S_l\!+\!\lambda_l E_{C_l}^{T}E_{C_l}+\bar{C}_l^{T}\left(W_{l}^{C}\!-\! W_{l}^{C}\left(\Psi_{l,k}+W_{l}^{C}\right)^{-1}W_{l}^{C}\right)\bar{C}_l  
        \\
        &\preceq S_l+\lambda_l E_{C_l}^{T}E_{C_l} + \bar{C}_l^{T}W_{l}^{C}\bar{C}_l= \bar{\Omega}_{l}.       
    \end{split}
\end{equation*}
Bounding $\Omega_{l,k}$ by $\bar{\Omega}_{l}$ and substituting in \eqref{eq::robust_riccati},
\begin{equation*}
    \begin{split}
        &P_{l,k+1} \\
        &= Q_{l} + F^T \Omega_{l,k}F - F^T\Omega_{l,k}G\left( R_{l} + G^T \Omega_{l,k} G\right)^{-1}G^T \Omega_{l,k}F 
        \\
        & \preceq Q_{l} + F^T \Omega_{l,k}F \preceq Q_{l} + F^T \bar{\Omega}_{l}F,
    \end{split}
\end{equation*}
the upper bound $P_{l,k+1}\preceq Q_{l} + F^T \bar{\Omega}_{l}F$ is obtained. 

\textbf{(3)} From (1) and (2), monotonicity of the Riccati recursion is guaranteed, i.e., $ 0=P_{l,0} \preceq P_{l,1} \preceq \dots \preceq Q_{l} + F^T \bar{\Omega}_{l}F,$ for all modes $l \in \hat{\Theta}$,  and the  convergence holds \cite{bellman1997introduction}.
\subsection{Proof of Lemma \ref{lemma::upper_bound_Lyapu}}
From the definition of $\Psi_{l}^{\delta} = \sum_{d \in \hat{\Theta}} \left( q_{ld} + \alpha_{ld} \Delta_{ld} \right) P_{d}$,  an upper bound can be found at $|\Delta_{ld}|=1$, 
\begin{equation*}
    \Psi_{l}^{\delta} \preceq \sum_{d \in \hat{\Theta}} \left( q_{ld} + \alpha_{ld} \right) P_{d} = \Psi_{l}.
\end{equation*}
Then, the inequality  $C_{l}^{\delta^T} \Psi_{l}^{\delta}  C_{l}^\delta \preceq C_{l}^{\delta^T}\Psi_{l} C_{l}^\delta$holds .  Given the structure $C_{l}^{\delta}=\bar{C}_l + M_{C_l}\Delta_{C_l} E_{C_l}$, an upper bound on $ C_{l}^{\delta^T}\Psi_{l} C_{l}^\delta$ can be obtained by directly using Lemma \ref{lemma::maximization}, as 
\begin{equation*}
    C_{l}^{\delta^T} \Psi_{l}^{\delta}  C_{l}^\delta \preceq C_{l}^{\delta^T}\Psi_{l} C_{l}^\delta \preceq \bar{C}_{l}^T \Psi_{l}^{C} \bar{C}_l + \lambda_{l} E_{C_l}^T E_{C_l},
\end{equation*}
with $\Psi_{l}^{C}=\Psi_{l} + \Psi_{l} M_{C_l}\left( \lambda_{l} I - M_{C_l}^T \Psi_{l} M_{C_l} \right)^{-1} M_{C_l}^T \Psi_{l} $, for $\lambda > \norm{M_{C_l}^T \Psi_{l} M_{C_l}}$. Multiplying these inequalities by $F_{CL}^T$ and $F_{CL}$ and subtracting $P_{l}$ from both sides, the result follows.

\subsection{Proof of Theorem \ref{thm::stability}}
In this proof, it is shown that the convergent Riccati solution of \eqref{eq::robust_riccati} satisfies the Lyapunov upper bound established in Lemma \ref{lemma::upper_bound_Lyapu}. Let the Riccati equation rewritten as
\begin{equation}\label{eq::Riccati_gain}
    F_{CL}^T \Omega_{l} F_{CL} - P_{l} = -\left( Q_{l} + K_{l}^T R_{l} K_{l} \right).
\end{equation}
For $\lambda_l>\norm{M_{C_l}^T \Psi_{l}M_{C_l}}$, the upper bound \eqref{eq::upper_bound_Lyapunov} is subtracted  from  \eqref{eq::Riccati_gain} which  ensures that the Riccati equation satisfies the Lyapunov upper bound, 
\begin{equation*}
    \begin{split}
        &\left( F_{CL}^T \Omega_{l} F_{CL} - P_{l} \right)-
        \\
        &\left(F_{CL}^T \bar{C}_{l}^T \Psi_{l}^{C} \bar{C}_{l}F_{CL} + \lambda_l F_{CL}^TE_{C_l}^T E_{C_l}F_{CL} - P_l \right) \succeq 0.
    \end{split}
\end{equation*}
Substituting $\Omega_{l}$ from \eqref{eq::aux_functions}, the inequality 

\begin{equation*}F_{CL}^T\left[ S_l+\bar{C}_l^{T}\left(\tilde{W}^{C}_{l} - \Psi_{l}^{C}\right) \bar{C}_{l} \right]F_{CL} \succeq 0,
\end{equation*}
holds, where $\tilde{W}^{C}_{l} = W_{l}^{C}-W_{l}^{C}\left(\Psi_{l}+W_{l}^{C}\right)^{-1}W_{l}^{C}$. It is therefore sufficient to impose $S_{l} \succeq \bar{C}_l^{T}\left(\Psi_{l}^{C} - \tilde{W}^{C}_{l}\right) \bar{C}_{l}$ to ensure that the Riccati solution dominates the derived upper bound. Since $Q_{l}, R_{l} \succ 0$, the difference $V_{k+1}-V_{k}$ satisfies the following inequalities
\begin{align*}
       x_{k}^T \left( F_{CL}^T C_{l}^{\delta^T} \Psi_{l}^{\delta}  C_{l}^\delta F_{CL}\!-\!P_{l} \right) x_{k} \leq-x_{k}^T\left( Q_{l} \!+\!K_{l}^T R_{l} K_{l} \right)x_{k} \!<\! 0,
\end{align*}
 except for $x_k = 0$, for all $l\in \hat{\Theta}$. From \cite[Theorem 3.9]{costa2005discrete}, the closed-loop system is MSS. 
\normalsize

\balance
\bibliographystyle{IEEEtran}
\bibliography{bib}

\end{document}